\magnification \magstep 1 
\input amssym.def     

\outer\def\beginsection#1\par{\vskip0pt plus.1\vsize\penalty-250
  \vskip0pt plus-.1\vsize\bigskip\vskip\parskip
  \message{#1}\leftline{\bf#1}\nobreak\smallskip\noindent}

\def\mapright#1{\smash{
   \mathop{\longrightarrow}\limits^{#1}}}
\def\mapleft#1{\smash{
   \mathop{\longleftarrow}\limits^{#1}}}
\def\mapdown#1{\Big\downarrow
   \rlap{$\vcenter{\hbox{$\scriptstyle#1$}}$}}
\def\mapup#1{\Big\uparrow
   \rlap{$\vcenter{\hbox{$\scriptstyle#1$}}$}}

\def\From{From}
\def\hhrule#1#2{\kern-#1
   \hrule height#1 depth#2 \kern-#2 }
\def\hvrule#1#2{\kern-#1{\dimen0=#1
    \advance\dimen0 by#2\vrule width\dimen0}\kern-#2 }
\def\makeblankbox#1#2{\setbox0=\hbox{A}
\hbox{\lower\dp0\vbox{\hhrule{#1}{#2}%
   \kern -#1
   \hbox to \wd0{\hvrule{#1}{#2}%
     \raise\ht0\vbox to #1{}
     \lower\dp0\vtop to #1{}
     \hfil\hvrule{#2}{#1}}
  \kern-#1\hhrule{#2}{#1}}}}
\def\qed{\makeblankbox{0pt}{.3pt}}

\def\AA{{\cal A}}
\def\Aut{{\rm Aut}}
\def\CC{{\cal C}}
\def\cc{{\Bbb C}}
\def\CCb{{\overline\CC}}
\def\ch{{\rm Ch}}
\def\chbar{{\overline{\ch}}}
\def\char{{\rm char}}
\def\colim{{\hbox{\rm colim}}}
\def\Cor{{\hbox{\rm Cor}}}
\def\DD{{\cal D}}
\def\ddot{{\cdot}}
\def\eg{{{\it e}.{\it g}.}}
\def\End{{\rm End}}
\def\fp{{\Bbb F}_p}
\def\fq{{\Bbb F}_q}

\def\gl{{\rm GL}}
\def\Hom{{\rm Hom}}
\def\ie{{{\it i}.{\it e}.}}
\def\ii{{\iota}}  
\def\Im{{\rm Im}}
\def\Ind{{\rm Ind}}
\def\mod{|}
\def\reg{{\rm reg}}
\def\Res{{\rm Res}}
\def\rr{{\Bbb R}}
\def\spec{{\rm Spec}}
\def\zz{{\Bbb Z}}
\def\v #1#2{{X_{#1}(#2)}}
\def\vv #1#2{{X'_{#1}(#2)}}
\def\proof{\par\noindent{\it Proof.}\quad}
\def\definition{\smallskip\par\noindent{\bf Definition.}\quad}
\def\remark{\smallskip\par\noindent{\bf Remark.}\quad} 
\def\example{\smallskip\par\noindent{\bf Example.}\quad}
\def\pra{\par}

\def\ada{1}
\def\alp{2}
\def\ati{3}
\def\atmc{4}
\def\ben{5}
\def\beninv{6}
\def\jfc{7}
\def\dww{8}
\def\ev{9}
\def\evtra{10}
\def\evbk{11}
\def\hls{12}
\def\hkr{13}
\def\lam{14}
\def\boaz{15}
\def\qui{16}
\def\pnilp{17}
\def\rec{18}
\def\ser{19}
\def\serbk{20}
\def\tho{21}
\def\ven{22}
\def\wilk{23}

\font \smallfont=cmr8 at 8pt
\font \smallbold=cmbx8 at 8pt

\centerline{\bf The spectrum of the Chern subring}
\vskip .2 in
\centerline{David J. Green and Ian J. Leary}  
\vskip .2 in 
\rightline{\it Dedicated to Charles Thomas, on}
\rightline{\it the occasion of his 60th birthday}
\vskip .2 in

{{\narrower\narrower\smallskip\noindent
{\smallbold Abstract.} {\smallfont For certain subrings of the
mod-${\scriptstyle p}$-cohomology of a compact Lie group, we give a
description of the spectrum, analogous to Quillen's description of the
spectrum of the whole cohomology ring.  Subrings to which our theorem
applies include the Chern subring.
Corollaries include a characterization of those groups for which the
Chern subring is F-isomorphic to the cohomology ring.  
}
\smallskip}}
\bigskip

\beginsection 1. Introduction. 

Let $G$ be a compact Lie group (\eg, a finite group) and let $H^*(G)=
H^*(BG;\fp)$ be its mod-$p$ cohomology ring.  This ring is a finitely
generated graded-commutative $\fp$-algebra.  In [\qui], D. Quillen
studied this ring from the viewpoint of commutative algebra.  His
results may be stated in terms of the prime ideal spectrum 
of $H^*(G)$, but the cleanest statement concerns the variety, 
$\v Gk$, of algebra homomorphisms from $H^*(G)$ to an algebraically 
closed field $k$ of characteristic $p$.  The Chern subring, 
$\ch(G)\subseteq H^*(G)$, is the subring generated by Chern classes of
unitary representations of $G$.  We give a description of $\vv Gk$, the
variety of algebra homomorphisms from $\ch(G)$ to $k$, analogous to
Quillen's description of $\v Gk$.  As corollaries of this result, we
classify the minimal prime ideals of $\ch(G)$, and characterize those
groups $G$ for which the natural map from $\v Gk$ to $\vv  Gk$ is a
homeomorphism.  

In the case when $G=E$ is an elementary abelian $p$-group, \ie, a
direct product of copies of the cyclic group of order $p$, $\v Ek$ is
naturally isomorphic to $E\otimes k$, where $E$ is viewed as a vector
space over $\fp$ and the tensor product is taken over $\fp$.  For
general $G$, Quillen describes $\v Gk$ as the colimit of the functor 
$(-)\otimes k$ over a category $\AA = \AA(G)$ with objects the
elementary abelian subgroups of $G$, and morphisms those group
homomorphisms that are induced by conjugation in $G$.  
Our description of $\vv Gk$ is as the colimit of the 
same functor over a category $\AA'$.  This category has the same
objects as Quillen's category, but a morphism in $\AA'$ is a group
homomorphism that merely preserves conjugacy in $G$.  In other words, 
a group homomorphism $f \colon  E_1 \rightarrow E_2$ satisfies:  
$$\eqalign{f\in \AA &\quad\Leftrightarrow\quad \exists g\, \forall e,
f(e) = g^{-1}eg,\cr
f\in \AA' &\quad\Leftrightarrow\quad \forall e\, \exists g,
f(e) = g^{-1}eg.}$$
This theorem is a corollary of a more 
general colimit theorem, which says, roughly, that the variety for any
subring of $H^*(G)$ that is both \lq large' and \lq natural' may be
expressed as such a colimit.  Other corollaries of this theorem,
in the case when $G$ is finite, include a description of the variety
for the subring of $\ch(G)$ generated by Chern classes of 
representations realizable over any subfield of $\Bbb C$, and 
a slight variation on the usual proof of Quillen's theorem 
in which transfers of Chern classes are used instead of the 
Evens norm map.  

\beginsection Acknowledgements.  

This work was started while the first named author held a Max Kade
visiting fellowship at the University of Chicago, and the second named
author held a Leibniz fellowship at the Max Planck Institut f\"ur
Mathematik, Bonn.  The first named author gratefully acknowledges the
support of the Deutsche Forschungsgemeinschaft and the Max Kade
Foundation.  The second named author gratefully acknowledges the
support of the Leibniz Fellowship Programme.  

The authors thank Dave Benson, Jon Carlson, Hans-Werner Henn, 
David Kirby and Jeremy Rickard for their helpful comments on this work.  

\beginsection 2. Representations and the Chern subring.  

First we recall some facts concerning Chern classes [\ati,\tho].  Let
$U(n)$ be the group of $n\times n$ unitary matrices, and $T(n)$ its
subgroup of diagonal matrices.  Then $H^*(T(n))$ is a free polynomial
algebra $\fp [x_1,\ldots, x_n]$ on $n$ generators of degree two.
$H^*(U(n))$ is isomorphic to $\fp[c_1,\ldots,c_n]$, where $c_i$ has
degree $2i$.  The map from $H^*(U(n))$ to $H^*(T(n))$ is injective and
sends $c_i$ to the $i$th elementary symmetric function in the $x_i$.  

Let $G$ be a compact Lie group.  $G$ has faithful finite-dimensional
complex representations, and any finite-dimensional 
representation is equivalent to a unitary representation.  If $\rho \colon 
G \rightarrow U(n)$ is a unitary representation of $G$, write $\rho \colon 
BG \rightarrow BU(n)$ for the induced map of classifying spaces, whose
homotopy class depends only on the equivalence class of~$\rho$.  The
$i$th Chern class of $\rho$ is defined by $c_i(\rho)= \rho^*(c_i)\in
H^{2i}(G)$.  Define $c_\ddot(\rho)$, the total Chern class of $\rho$,
to be $1+ c_1(\rho) + \cdots + c_n(\rho)$.  Chern classes enjoy the
following properties (`Whitney sum formula' and `naturality'), for any
$\theta \colon  G \rightarrow U(m)$ and any $f \colon  H\rightarrow G$:  
$$c_\ddot(\rho\oplus\theta) = c_\ddot(\rho)c_\ddot(\theta) \qquad 
c_\ddot(\rho\circ f)= f^*c_\ddot(\rho).$$
There is a unique way to define Chern classes for virtual
representations so that they continue to enjoy the above properties.
Let $c_\ddot'= 1+ c'_1+ c'_2 + \cdots$ be the unique power series in
$\fp[[c_1,\ldots,c_n]]$ satisfying $c'_\ddot c_\ddot=1$, and define
$c_i(-\rho)= \rho^*(c'_i)$.  In general infinitely many of the
$c_i(-\rho)$ will be non-zero, but note that they are all expressible
in terms of the $c_i(\rho)$.  

\definition
The Chern subring $\ch(G)$ of $H^*(G)$ is the subring generated by the
$c_i(\rho)$ for all $i$ and all virtual representations $\rho$.  
\medskip

By the above remarks, the $c_i(\rho)$ as $\rho$ ranges over the
irreducible representations of $G$ suffice to generate $\ch(G)$.  In
the case when $G$ is finite, it follows that $\ch(G)$ is finitely
generated, since $G$ has only finitely many inequivalent irreducible
representations.  For general $G$ it is also true that $\ch(G)$ is
finitely generated.  This is a special case of the following
proposition.  

\proclaim Proposition 2.1.  
Let $G$ be a compact Lie group, and $\rho \colon  G \rightarrow U(n)$ a 
faithful unitary representation of $G$.  If $R$ is a subring of
$H^*(G)$ containing each $c_i(\rho)$, then $R$ is finitely generated.  

\proof 
Venkov showed that $H^*(G)$ is finitely generated by showing that 
$H^*(G)$ is finite over (\ie, is a finitely generated module for)
$H^*(U(n))$ 
[\ven,\qui,\ben].  Now $R$ is an $H^*(U(n))$-submodule of $H^*(G)$, so
is finitely generated since $H^*(U(n))$ is Noetherian.  
\qed

\remark For $G$ finite, the finite-generation of $H^*(G)$ is due
independently to Evens and to Venkov by completely different proofs 
[\ev,\ven].  There is another proof (closer to Evens' than to Venkov's) 
in [\dww].   

\definition
A virtual representation $\rho$ of $G$ is said to be $p$-regular if
the virtual dimension of $\rho$ is strictly positive and for every
elementary abelian subgroup $E\cong (\zz/p)^n$ of $G$, the restriction
to $E$ of $\rho$ is a direct sum of copies of the regular
representation of $E$.  

\proclaim Proposition 2.2. 
For each prime $p$, every compact Lie group $G$ has a 
$p$-regular representation.  

\proof 
$G$ has a faithful representation in $U(n)$ for some $n$, and every
elementary abelian subgroup of $U(n)$ is conjugate to a subgroup of 
$T(n)$, the torus consisting of diagonal matrices.  Thus it suffices
to show that $U(n)$ has a virtual representation whose restriction to 
$(\zz/p)^n\subseteq T(n)$ is the regular representation.  

Recall that 
the representation ring $R(T(n))$ of $T(n)$ is isomorphic to the
Laurent polynomial ring
$\zz[\tau_1,\tau_1^{-1},\ldots,\tau_n,\tau_n^{-1}]$, where $\tau_i$ is
the 1-dimensional representation 
$$\tau_i \colon  {\rm diag}(\xi_1,\ldots,\xi_n)\mapsto \xi_i.$$
$R(U(n))$ maps injectively to $R(T(n))$ with image the subring
$\zz[\sigma_1,\sigma_2,\ldots,\sigma_n,\sigma_n^{-1}]$, where
$\sigma_i$ is the $i$th elementary symmetric function in the
$\tau_j$.  

The polynomial 
$$P = \prod_{i=1}^n (1 + \tau_i + \cdots + \tau_i^{p-1})$$
is a symmetric polynomial in the $\tau_j$, and so is expressible in
terms of $\sigma_1,\ldots,\sigma_n$.  The corresponding ($p^n$-dimensional)
representation of $U(n)$ restricts to $(\zz/p)^n\subseteq T(n)$ as the
regular representation.   
\qed

\remark 
For $G$ finite, the regular representation of $G$ is of course
$p$-regular.  

Using Quillen's result that we state as Theorem~5.1, it may be shown
that $H^*(G)$ is finite over the subring generated by the Chern
classes of any $p$-regular representation.  For genuine (as opposed to
virtual) representations, this can be deduced from Venkov's result:  
If $\rho \colon  G \rightarrow U(n)$ is a $p$-regular 
representation of $G$, the kernel of $\rho$ contains no elements of
order $p$, and is therefore a finite group of order coprime to $p$.  
It follows that $H^*(\rho(G))\cong H^*(G)$ (consider the spectral
sequence for the extension $\ker(\rho)\rightarrow G \rightarrow
\rho(G)$), and hence $H^*(G)$ is finite over the image of $\rho^*$.   
When $p=2$, the representation constructed in the 
proof of Proposition~2.2 is a genuine representation.  David Kirby has
shown us an argument to prove that $U(n)$ has a $p$-regular
genuine representation if and only if either $p=2$, or $n=1$, or
$(p,n)=(3,2)$.  

\beginsection 3. Varieties for cohomology. 

Let $k$ be an algebraically closed field of characteristic $p$, and
let $R$ be a finitely generated commutative $\fp$-algebra.  Define
$V_R(k)$, the variety for $R$, to be the set of ring homomorphisms 
from $R$ to $k$, with the Zariski topology, \ie, the smallest
topology in which the set
$$F_I=\{\phi\colon  R\rightarrow k \mid \ker(\phi)\supseteq I\}$$
is closed for each ideal $I$ of $R$.  A ring homomorphism $f \colon  R
\rightarrow S$ gives rise to a continuous map $f^*\colon  V_S(k) \rightarrow
V_R(k)$.  If $S$ is finite over $f(R)$ (\ie, $S$ is a finitely generated
$f(R)$-module) then $f^*$ is a closed mapping and has finite fibres,
by the `going up' or `lying over' theorem [\atmc,\beninv].  If $S$ is
finite over $f(R)$ and $\ker(f)$ is nilpotent, then $f^*$ is
surjective.  

There is a continuous map from $V_R(k)$ to $\spec(R)$, the prime ideal
spectrum of $R$, that sends the map $\phi \colon  R\rightarrow k$ to the
ideal $\ker(\phi)$.  If the transcendence degree of $k$ over $\fp$ is
sufficiently large (as large as a generating set for $R$ will
suffice), then this map is surjective.  Thus information about
$V_R(k)$ gives rise to information about $\spec(R)$.  

We shall also require the following:  
\proclaim Proposition 3.1.  
(a) Let $S$ be a subring of $R$
containing $R^p$, the subring of $p$th powers of elements of $R$.
Then the natural map from $V_R(k)$ to $V_S(k)$ is a homeomorphism.   
\pra 
\noindent 
(b) Let a finite group $G$ act on $R$, with fixed point subring
$S=R^G$.  Then the natural map $V_R(k)\rightarrow V_S(k)$ induces a
homeomorphism $V_R(k)/G \rightarrow V_S(k)$.  

\proof
Each of these claims 
may be proved by showing that $R$ is finite over $S$, and
deducing that the map given is continuous, closed, and a
bijection.  See for example [\atmc,\beninv,\qui].  
\qed

For $p=2$, the ring $H^*(G)$ is commutative.  For $p$ an odd prime,
elements of $H^*(G)$ of odd degree are nilpotent, so although $H^*(G)$
is not commutative, the quotient of $H^*(G)$ by its radical, $h^*(G)=
H^*(G)/\sqrt 0$, is commutative.  Any homomorphism from $H^*(G)$ to
$k$ factors through $h^*(G)$.  Define $X_G(k)$ to be the variety
$V_R(k)$ for $R= h^*(G)$.  By the above remark, points of $X_G(k)$ may
be viewed as homomorphisms from $H^*(G)$ to $k$.  Let $S$ be the
subring of elements of $H^*(G)$ of even degree.  A homomorphism
$\phi\colon S\rightarrow k$ extends uniquely to a homomorphism from $H^*(G)$
to $k$ (if $x$ is in odd degree, then either $p$ is odd, and $x^2=0$, 
or $p=\char(k)=2$, so in either case $\phi(x)$ is the unique square
root of $\phi(x^2)$).  It follows that the natural map
$X_G(k)\rightarrow V_S(k)$ is a homeomorphism, since it is a closed,
continuous bijection.  Hence $X_G(k)$ could equally be defined in
terms of $S$.  

Note that a group homomorphism $f\colon  H\rightarrow G$ induces a map $f_*\colon 
X_H(k) \rightarrow X_G(k)$.  We write $\ii^G_H$ for $f_*$ in the case
when $f$ is the inclusion of a subgroup $H$ in $G$.  A theorem of
Evens and Venkov [\ev,\ven] states that in this case $H^*(H)$ is
finite over $H^*(G)$.  (To deduce this from the result quoted in the
proof of Proposition~2.1, note that a faithful representation of $G$
restricts to a faithful representation of $H$.)  It follows that
$\ii^G_H$ is closed and has finite fibres.  

Define $X'_G(k)$ to be $V_{\ch(G)}(k)$.  By Proposition~2.1 the
natural map from $X_G(k)$ to $X'_G(k)$ is surjective, closed, and has
finite fibres.  

\proclaim Proposition 3.2.  Let $\rho$ be a representation of 
$G$, and let $R(n)$ be the subring of $H^*(G)$ generated by the Chern
classes of $n\rho= \rho^{\oplus n} = 
\rho \oplus\cdots\oplus \rho$.  Then the natural
map $V_{R(1)}(k) \rightarrow V_{R(n)}(k)$ is a homeomorphism.  

\proof If $p$ does not divide $n$, then 
$$c_i(n\rho) = nc_i(\rho) + P(i,n),$$ 
for some expression $P(i,n)$ in the $c_j(\rho)$ for $j<i$.  So in this
case $R(n) = R(1)$.  On the other hand, if $n=pm$ then 
$$c_\ddot(n\rho) = c_\ddot(m\rho)^p = 1 + c_1(m\rho)^p + c_2(m\rho)^p
+ \cdots,$$ 
so in this case $R(n) = R(m)^p$, the subring of $p$th powers of
elements of $R(m)$, and the map $V_R(k)\rightarrow V_{R^p}(k)$ is a
homeomorphism.  
\qed

The methods that we shall use to study the Chern subring apply equally
to the Stiefel-Whitney subring, defined analogously, in the case when
$p=2$.  (For information concerning Stiefel-Whitney classes see
[\tho]).  As an alternative, the following proposition may be
applied.  

\proclaim Proposition 3.3. Let $p=2$ and let $S$ be the subring of
$H^*(G)$ generated by Stiefel-Whitney classes of real representations
of $G$.  Then 
$$S^2\subseteq \ch(G)\subseteq S,$$
and the natural map from $V_S(k)$ to $X'_G(k)$ is a homeomorphism.  

\proof  If $\theta$ is an $n$-dimensional real representation of $G$,
then $\theta^\cc$, the complexification of $\theta$, is an
$n$-dimensional complex representation of $G$ with $c_i(\theta^\cc)=
w_i^2(\theta)$.  Conversely, if $\psi$ is an $n$-dimensional complex
representation of $G$ and $\psi_\rr$ is the same representation viewed
as a $2n$-dimensional real representation, then $w_i(\psi_\rr)=0$ for
$i$ odd and $w_{2i}(\psi_\rr)= c_i(\psi)$.  This proves the claimed
inclusions.  The claimed homeomorphism follows from
Proposition~3.1(a).  
\qed

\beginsection 4. Examples.  

In this section we discuss the case of an elementary abelian
$p$-group, and also give an example to show that the map 
$X_G(k)\rightarrow X'_G(k)$ is not always a homeomorphism.  
This example was the starting point for the work of this paper.  

Let $E$ be an elementary abelian $p$-group of rank $n$, $E\cong
(\zz/p)^n$.  Then $E$ may be viewed as a vector space over $\fp$.
Write $E^*$ for $\Hom(E,\fp)$.  There is a natural isomorphism
$E^*\cong H^1(E)$.  For $p=2$, $H^*(E)$ is isomorphic to the
symmetric algebra on $H^1(E)$, or equivalently, the ring of polynomial
functions on $E$ viewed as a vector space:  
$$H^*(E)\cong S(E^*)\cong \fp[E].$$ 
For $p>2$, the Bockstein $\beta \colon  H^1(E)\rightarrow H^2(E)$ is
injective, and $H^*(E)$ is isomorphic to the tensor product of the
exterior algebra on $H^1(E)$ tensored with the symmetric algebra on $B
= \beta(H^1(E))$:  
$$H^*(E)\cong \Lambda(E^*)\otimes S(B)\cong \Lambda(E^*)\otimes
\fp[E].$$
In any case, $h^*(E)$ is naturally isomorphic to $\fp[E]$, generated
in degree one for $p=2 $ and in degree two for odd $p$.  It follows
that $X_E(k)$ is naturally isomorphic to $E\otimes k$, where $E$ is
viewed as a vector space over $\fp$ and the tensor product is over
$\fp$, so that $E\otimes k \cong k^n$.  

Irreducible representations of $E$ are 1-dimensional, and the map
$\rho \mapsto c_1(\rho)$ is a natural bijection between the set of
irreducible representations of $E$ and $B=\beta(H^1(E))$.  (When
$p=2$, $\beta = {\rm Sq}^1$, and so $\beta(x)= x^2$.)  The Chern
subring $\ch(E)$ of $H^*(E)$ is the subalgebra of $H^*(E)$ generated
by $B$.  For $p>2$ this subring maps onto $h^*(E)$, and for $p=2$ it
maps onto $h^*(E)^2$, the subring of squares of elements of $h^*(E)$.
In any case, the map from $X_E(k)$ to $X'_E(k)$ is a homeomorphism.  

\proclaim Proposition 4.1.  Let $\rho$ be a direct sum of copies of the
regular representation of $E$, and let $R$ be the subring of $H^*(E)$
generated by the Chern classes of $\rho$.  Then the natural map from
$X_E(k)$ to $V_R(k)$ factors through a homeomorphism 
$$k^n/\gl_n(\fp)\cong X_E(k)/\gl(E) \rightarrow V_R(k).$$

\proof By Proposition 3.2 it suffices to consider the case when $\rho$
is the regular representation.  Identify $\ch(E)$ with $\fp [E]$,
generated in degree two.  The total Chern class of $\rho$ is 
$$c_\ddot(\rho) = \prod_{x\in E^*}(1+x).$$ 
This is invariant under the full automorphism group, $\gl(E)$, of
$E$.  By a theorem of Dickson, the only $i>0$ for which $c_i(\rho)$ is
non-zero are $i=p^n-p^j$, where $0\leq j<n=\dim_{\fp}(E)$.  Moreover,
these $c_i(\rho)$ freely generate a polynomial subring of $\fp[E]$,
and this is the complete ring of $\gl(E)$-invariants in $\fp[E]$
[\beninv,\wilk].  The claim follows by part (b) of Proposition~3.1.  
\qed

\remark Let $A$ be a non-identity element of $\gl_n(\fp)$, and let $v$
be an element of $k^n$ fixed by $A$.  Then $v$ is in the kernel of
$I-A$, a non-zero matrix with entries in $\fp$, and so $v$ lies in a
proper subspace of $k^n$ defined over $\fp$ (\ie, a subspace of the
form $V\otimes k$ for some proper $\fp$-subspace $V$ of $\fp^n$).  It
follows that $\gl_n(\fp)$ acts freely on the complement of all such
subspaces.  For an elementary abelian group $E$, let 
$$X^+_E(k) = X_E(k)\setminus \bigcup_{F<E}\ii^E_F(X_F(k)).$$
By the above argument, $\gl(E)$ acts freely on $X^+_E(k)$.

\example Let $q=p^n$ for some $n\geq 2$.  Let $G$ be the affine
transformation group of the line over $\fq$.  Then $G$ is expressible
as an extension with kernel $E=(\fq,+)$, an elementary abelian
$p$-group of rank $n$, and quotient $Q = \gl_1(\fq)$, cyclic of order
$q-1$.  The conjugation action of $Q$ on $E$ is transitive on
non-identity elements of $E$.  One example of such a group is the 
alternating group $A_4$ ($p=2$, $n=2$).  
An easy transfer argument shows that
$H^*(G)$ maps isomorphically to the ring of invariants, $H^*(E)^Q$,
and it follows from Proposition~3.1 that $X_G(k)$ is homeomorphic to
$X_E(k)/Q = k^n/Q$, where $Q= \gl_1(\fq)\leq \gl_n(\fp)\leq
\gl_n(k)$.  

It is easy to check that $G$ has exactly $q$ distinct irreducible
representations.  All but one of these are 1-dimensional and restrict
to $E$ as the trivial representation.  The other one is
$(q-1)$-dimensional and restricts to $E$ as the regular representation
minus the trivial representation.  Hence by Proposition~4.1, $X'_G(k)$
is homeomorphic to $X_E(k)/\gl(E)=k^n/\gl_n(\fp)$.  Thus the map from
$X_G(k)$ to $X'_G(k)$ is not a homeomorphism.  

\beginsection 5. Quillen's colimit theorem.  

In [\qui], Quillen showed that for general $G$, $X_G(k)$ is determined
by the elementary abelian subgroups of $G$.  Roughly speaking, he
showed that $X_G(k)$ is equal to the union of the images of the
$X_E(k)$, where $E$ ranges over the elementary abelian subgroups of
$G$, and that as little identification takes place between the points
of the $X_E(k)$ as is consistent with the fact that inner
automorphisms of $G$ act trivially on $H^*(G)$.  More precisely, 
let $f\colon E_1 \rightarrow E_2$ be a homomorphism between elementary
abelian subgroups of $G$ that is induced by an inner automorphism of
$G$.  Then the following diagram commutes.  
$$\matrix{H^*(G)&\mapleft{{\rm Id}}&H^*(G)\cr 
\mapdown{\Res}&&\mapdown{\Res}\cr
H^*(E_1)&\mapleft{f^*}&H^*(E_2)\cr}$$
Consequently the following diagram commutes.  
$$\matrix{X_G(k)&\mapright{{\rm Id}}&X_G(k)\cr 
\mapup{\ii}&&\mapup{\ii}\cr
X_{E_1}(k)&\mapright{f_*}&X_{E_2}(k)\cr}$$
This fact motivates the following definition.  

\definition The Quillen category $\AA$ for a compact Lie group $G$ 
and a prime $p$ is 
the category whose objects are the elementary abelian $p$-subgroups of
$G$, with morphisms from $E_1$ to $E_2$ being those group
homomorphisms that are induced by conjugation in $G$.  
Any such group homomorphism is of course injective.  
\medskip

In general $G$ will have infinitely many elementary abelian
$p$-subgroups.  These subgroups form finitely many conjugacy classes
though ([\qui], lemma~6.3).  Thus although the Quillen category for
$G$ is infinite in general, it contains only finitely many isomorphism
types of object (or is \lq skeletally finite').  

The morphisms $f\colon E_1\rightarrow E_2$ in the Quillen category are
precisely the maps for which the diagram above commutes.  It follows
that the natural map
$$\coprod_{E\leq G \atop E\,\,{\rm el.\,ab.}} X_E(k)\longrightarrow
X_G(k)$$  
factors through a map $\alpha\colon  \colim_\AA X_E(k)\rightarrow X_G(k)$.  

\proclaim Theorem 5.1.  (Quillen [\qui]) The map $\alpha \colon \colim_\AA
X_E(k)\rightarrow X_G(k)$ is a homeomorphism.  

The map $\alpha$ is continuous and is closed because  $\AA$ is 
skeletally finite.  Thus the main content of the theorem is that
$\alpha$ is a bijection.  We shall use only half of this theorem, 
the statement that $\alpha$ is surjective, in our main theorem.  
The surjectivity of $\alpha$ is equivalent to the statement \lq an
element of $H^*(G)$ is nilpotent if and only if its image in each
$H^*(E)$ is nilpotent'.  

\beginsection 6. A new colimit theorem.  

Motivated by the Quillen category, we define: 

\definition A category of elementary abelian subgroups of $G$ is a
category whose objects are (all of) the elementary abelian $p$-subgroups 
of $G$, and whose morphisms from $E_1$ to $E_2$ are injective group 
homomorphisms.  
\medskip

The Quillen category, $\AA(G)$, is of course a category of elementary
abelian subgroups of $G$.  Another example is the category
$\CC_\reg(G)$, 
with the morphism set $\CC_\reg(E_1,E_2)$ equal to the
set of all 1-1 group homomorphisms from $E_1$ to $E_2$.  Any category
of elementary abelian subgroups of $G$ is a subcategory of
$\CC_\reg(G)$.  

Any subring $R$ of $H^*(G)$ gives rise to a category $\CC(R)$ of
elementary abelian subgroups of $G$, where $f\colon E_1\rightarrow E_2$ is a
morphism in $\CC(R)$ if and only if 
$$\matrix{R&\mapleft{{\rm Id}}&R\cr 
\mapdown{\Res}&&\mapdown{\Res}\cr
h^*(E_1)&\mapleft{f^*}&h^*(E_2)\cr}$$
commutes.  Note that we use $h^*(E_i)$, the cohomology ring modulo its
radical, rather than $H^*(E_i)$.  Note that if $f$ is an
isomorphism of groups and is a morphism in $\CC(R)$, $f^{-1}$ is also
in $\CC(R)$.  Each $\CC(R)$ contains the Quillen category, and hence
is skeletally finite.  By the argument given in Section~5, the map 
$\alpha\colon \colim_\AA X_E(k)\rightarrow X_G(k)$ induces a map
$\gamma=\gamma(R)\colon  \colim_{\CC(R)}X_E(k)\rightarrow V_R(k)$.  

\definition Say that a subring of $H^*(G)$ is large if it contains the
Chern classes of some $p$-regular representation of $G$.  Say that a
subring of $H^*(G)$ is natural if it is generated by homogeneous
elements and is closed under the action of the Steenrod algebra.  
\medskip

The new colimit theorem of the title of this section is:  

\proclaim Theorem 6.1.  Let $G$ be a compact Lie group, and let $R$ be
a subring of $H^*(G)$ that is both large and natural.  Then the map 
$$\gamma\colon \colim_{\CC(R)}X_E(k)\rightarrow V_R(k)$$ 
is a homeomorphism.  

It is possible that this theorem could be proved using more general 
colimit theorems due to S.~P. Lam, to D. Rector, and to 
H.-W. Henn, J. Lannes and L. Schwartz [\lam,\rec,\hls].  These
theorems say, roughly speaking, that the variety for any Noetherian
algebra over the Steenrod algebra should be expressible as a similar
sort of colimit.  Even with these theorems, 
Quillen's description of $X_G(k)$ would still be 
needed to identify the categories that arise 
with categories of elementary abelian subgroups of $G$.  The proof 
given below is more elementary, 
in that it relies on no work that is more recent than that of
Quillen.  

The proof of the theorem uses the following lemma.  

\proclaim Lemma 6.2.  Let $S$ be the subring of $H^*(G)$ generated by
the Chern classes of some $p$-regular representation of $G$.  Then
$\CC(S)$ is equal to the category $\CC_\reg$ defined above, and the
map $\gamma(S)$ is a homeomorphism.  

\proof Let $F$ be a maximal elementary abelian subgroup of $G$.  
Note that the natural map from $X_F(k)$ (mapping the category with
one object and one morphism to $\CC_\reg$) induces a homeomorphism 
$$X_F(k)/\gl(F) \cong \colim_{\CC_\reg} X_E(k).$$
By Proposition~4.1, the image of $\ii^G_F \colon  X_F(k)\rightarrow V_S(k)$ 
is homeomorphic to $X_F(k)/\gl(F)$.  If $E$ is any elementary abelian
subgroup of $G$ 
and $f \colon E \hookrightarrow F$ is any injective group homomorphism,
then $\Res^G_E(\rho) $ and $f^*\Res^G_F(\rho)$ are equal to 
a sum of (the same number of) copies of the regular representation of
$E$.  Hence $\CC(S)$ is equal to $\CC_\reg$, and $\Im(\gamma)=
\Im(\ii^G_F)$.  It follows that $\gamma$ is a homeomorphism onto its
image.  Finally, by Theorem~5.1, this image is the whole of $V_S(k)$.
\qed

\par\noindent 
{\it Proof of the theorem.}\quad Since $R$ is large, it contains a
subring $S$ satisfying the conditions of Lemma~6.2.  Let $E$, $F$ be
two elementary abelian subgroups of $G$, suppose that the rank of $E$
is less than or equal to that of $F$, and suppose that 
$\phi\in X_E(k)$ and $\psi\in X_F(k)$ define the same point of
$V_R(k)$.  {\it A fortiori\/} 
$\phi$ and $\psi$ define the same point of $V_S(k)$,
and so by Lemma~6.2 there is an injective group homomorphism 
$f\colon  E \hookrightarrow F$ such that $\psi = \phi\circ f^* =
f_*(\phi)$.  It suffices to show that such an $f$ is in $\CC(R)$.  

For any such $f \colon  E \hookrightarrow F$, let $\cal S$ be the set 
of subgroups of $E$ such that $f$ restricted to $E$ is a morphism in
$\CC(R)$: 
$${\cal S} = \{ E'\leq E \mid \left(f\mod_{E'}\colon E'\rightarrow
F\right) \in \CC(R)\},$$ 
and define a subset $X(f)$ of $X_E(k)$ by 
$$X(f) = \{ \phi \in X_E(k) \mid f_*(\phi)\circ\Res^G_{E'}\mod_R 
= \psi\circ\Res^G_F\mod_R \}. $$
\From\ the definitions, 
$$X(f) \supseteq \bigcup_{E' \in \cal S} \ii^E_{E'} X_{E'}(k),$$ 
and it suffices to show that equality holds.  Note that 
a subgroup $E'\leq E$ is in $\cal S$ if and only
if $\ii X_{E'}(k)$ is a subset of $X(f)$.  Hence it suffices to show
that $X(f)$ is equal to some union of sets of the form $\ii X_{E'}(k)$.  
Now let $I(f)$ be the ideal of $H^*(E)$ generated by all elements of
the form $\Res^G_E(r)- f^*\Res^G_F(r)$, where $r\in R$.  The
subvariety of $X_E(k)$ defined by $I(f)$ is the set $X(f)$ defined
above.  Since $R$ is natural (in the sense defined above the
statement), the ideal $I(f)$ is homogeneous and closed under the
action of the Steenrod algebra.  But by a theorem of Serre
[\ser,\qui], the variety corresponding to any such ideal of $H^*(E)$
has the required form.  
\qed

Minimal prime ideals of a commutative $\fp$-algebra $R$ correspond to
irreducible components of $V_R(k)$.  Hence one obtains:  

\proclaim Corollary 6.3.  Let $R$ be a large, natural subring of
$H^*(G)$.  The minimal prime ideals of $R$ are in bijective
correspondence with the isomorphism types of maximal objects in
$\CC(R)$.  

An object of a category is called maximal if every map 
from it is an isomorphism.  An
isomorphism class of maximal objects in the Quillen category is a
conjugacy class of maximal elementary abelian subgroups of $G$.  

\proclaim Corollary 6.4.  Let $R$ and $S$ be large natural subrings of
$H^*(G)$, and suppose that $R$ is a subring of $S$.  The natural map
$V_S(k)\rightarrow V_R(k)$ is a homeomorphism if and only if the
categories $\CC(R)$ and $\CC(S)$ are equal.  

\proof A direct proof could be given at this stage, but it is easier
to apply Proposition~9.1, which implies that no subcategory of
$\CC_\reg$ that strictly contains $\CC(S)$ gives rise to the same
variety.  \qed

\beginsection 7.  Applications to rings of Chern classes.  

We start by defining some categories of elementary abelian subgroups
of $G$.  

\definition Define categories of elementary abelian subgroups of $G$, 
$\AA'$, $\AA'_\rr$, $\AA'_{\rm P}$, and for each $d$ dividing $p-1$, 
$\AA_d$, 
by stipulating that an 
injective homomorphism $f \colon  E \rightarrow F$ is in:  
\pra $\AA'$ if $\forall e$, $f(e)$ is conjugate (in $G$) to $e$; 
\pra $\AA'_\rr $ if $\forall e$, $f(e)$ is conjugate to $e$ or to
$e^{-1}$; 
\pra $\AA'_{\rm P}$ if $\forall e$, the subgroups $\langle e \rangle $
and $\langle f(e) \rangle$ are conjugate; 
\pra $\AA'_d$ if $\forall e$, $f(e)$ is conjugate to $\xi(e)$ for some
$\xi$ in the order $d$ subgroup of $\Aut(\langle e\rangle)$.  
\medskip

Note that for $p=2$, $\AA'= \AA'_\rr = \AA'_{\rm P}$, and for odd $p$, 
$\AA'_\rr = \AA'_2$, and $\AA'_P = \AA'_{p-1}$.  Note also that the
difference between $\AA'$ and the Quillen category $\AA$ is the
difference between \lq $\forall e\, \exists g f(e)= g^{-1}eg$'
and \lq $\exists g\, \forall e f(e)= g^{-1}eg$'.  The reason for 
introducing these categories is the following proposition.  

\proclaim Proposition 7.1.  Let $K$ be a subfield of $\cc$ and let 
$|K(\zeta_p)\colon K|=d$, where $\zeta_p$ is a primitive $p$th root of 1.  
Let $G$ be a compact Lie group, and in
cases (c) and (d) suppose that $G$ is finite.  Let $R$ be the subring
of $H^*(G)$ generated by Chern classes of:  
\pra 
(a) All representations of $G$; 
\pra 
(b) Representations of $G$ realisable over the reals; 
\pra 
(c) Permutation representations of $G$;  
\pra 
(d) Representations of $G$ realisable over $K$.  
\pra 
In each case, the variety $V_R(k)$ is homeomorphic to $\colim_{\CC(R)}
X_E(k)$.  The category $\CC(R)$ is:  
$$\hbox{(a) $\AA'$,}\quad 
\hbox{(b) $\AA'_\rr$,}\quad
\hbox{(c) $\AA'_{\rm P}$,}\quad
\hbox{(d) $\AA'_d$.}$$   

\proof In each case, the morphisms in the category given are precisely
those group homomorphisms for which $\chi(e) = \chi(f(e))$ for all 
characters $\chi$ coming from representations of the given type.  
(See [\serbk] Chapter~12 for case (d), and for case (c) note that if 
$e$, $e'$ are elements of $G$ of order $p$, then the permutation 
action of $e$ on $G/\langle e'\rangle$ has a fixed point if and only
if $\langle e\rangle$ is conjugate to $\langle e'\rangle$.)  The 
proposition therefore follows from the lemma below.  \qed 

\proclaim Lemma 7.2.  
Let $A$ be an additive subgroup of the representation ring of $G$, 
generated by genuine representations, and containing a $p$-regular
representation.  Let $R$ be the subring of $H^*(G)$ generated by the
Chern classes of all elements of $A$.  Then $R$ is large and natural, 
and hence by Theorem~6.1 
$$\gamma \colon  \colim_{\CC(R)} X_E(k) \rightarrow V_R(k)$$ 
is a homeomorphism.  Furthermore, $f \colon  E \hookrightarrow F$ is a
morphism in $\CC(R)$ if and only if for all $e\in E$, and all 
characters $\chi$ of elements of $A$, $\chi(e) = \chi(f(e))$.  

\proof First, suppose that $A$ is generated by a single representation
$\rho$.  The image of $\rho^* \colon  H^*(U(n)) \rightarrow H^*(G)$ is
natural since $H^*(U(n))$ is graded and acted upon by the Steenrod
algebra.  The general case follows from the Cartan formula.  By
hypothesis, $R$ is large.  The claimed homeomorphism follows from
Theorem~6.1, and it only remains to describe $\CC(R)$.  

A representation is determined up to equivalence by its character.  
Hence for any $f$ as in the statement and $\rho$ a generator of $A$, 
$f^*\Res^G_F c_\ddot(\rho)= \Res^G_E c_\ddot(\rho)$, and so any such
$f$ is in $\CC(R)$.  For the converse, note that since $\ch(E)$ is a 
unique factorization domain, a representation of $E$ is determined up
to equivalence by its dimension and its total Chern class.  Thus if 
$f \colon  E \hookrightarrow F$ is a homomorphism for which there exists $e$
and $\chi$ with $\chi(f(e))\neq \chi(e)$, there 
exists $i$ and $\rho$, a generator of $A$, such that
$\Res^G_E(c_i(\rho)) - f^*\Res^G_F(c_i(\rho))\neq 0$.  Hence $f$ is
not in $\CC(R)$.  
\qed


Quillen's description of $X_G(k)$ (Theorem 5.1, and Theorem~8.1),
Corollary~6.4 and Proposition~7.1 together yield: 

\proclaim Corollary 7.3.  The natural map $X_G(k)\rightarrow X'_G(k)$
is a homeomorphism if and only if the categories $\AA(G)$ and
$\AA'(G)$ are equal.  \qed

\example (A $p$-group $G$ for which the map $X_G(k)\rightarrow
X'_G(k)$ is not a  homeomorphism.)  
Let $E$ be the additive group $\fp^n$ for some $n>2$, 
and let $U\leq \gl(E)= \gl_n(\fp)$
be the Sylow $p$-subgroup of $\gl(E)$ consisting of upper triangular
matrices.  Let $Q$ be the subgroup of $U$ consisting of all matrices
$(a_{i,j})$ that are constant along diagonals, \ie, $a_{i,j} =
a_{i+1,j+1}$ whenever $1\leq i < m$ and $1\leq j <m$.  Finally, let
$G$ be the split extension with kernel $E$ and quotient $Q$.  
The group $E$ is a maximal elementary abelian subgroup of $G$.  Easy
matrix calculations show that the orbits of the action of $Q$ on
elements of $E$ are equal to the orbits of the action of $U$, and 
that any element of $\gl(E)$ that preserves the $U$-orbits in $E$ is
fact an element of $U$.  It follows that the image of $X_E(k)$ in 
$X'_G(k)$ is $X_E(k)/U$, whereas the image of $X_E(k)$ in $X_G(k)$ is
of course $X_E(k)/Q$.  Thus $G$ is a $p$-group such that the fibres
of the map $X_G(k)\rightarrow X'_G(k)$ above general points of one 
irreducible component have order $|U\colon Q|= p^{(n-1)(n-2)/2}$.  

\example (A group $G$ for which $X'_G(k)$ has fewer irreducible
components than $X_G(k)$.)  Let $G$ be $\gl_3(\fp)$.  There are two
distinct Jordan forms for elements of order $p$ in $G$ (resp.\ one if 
$p=2$), and hence $G$ has two conjugacy classes (resp.\ one conjugacy
class if $p=2$) of elements of order $p$.  All maximal elementary abelian
subgroups of $G$ have rank two.  The subgroups 
$$E_1=\pmatrix{1&*&*\cr 0&1&0\cr 0&0&1\cr}\quad 
\hbox{and}\quad E_2=\pmatrix{1&0&*\cr 0&1&*\cr 0&0&1\cr}$$
are maximal elementary abelian subgroups, and are not conjugate,
although every non-identity element of $E_1$ is conjugate to every
non-identity element of $E_2$.  It follows that the images of
$X_{E_1}(k)$ and $X_{E_2}(k)$ in $X_G(k)$ are distinct irreducible
components of $X_G(k)$, whereas their images in $X'_G(k)$ are equal.  

\beginsection 8.  Transfers of Chern classes.  

Throughout this section, $G$ shall be a finite group.  
Following Moselle [\boaz], we consider the \lq Mackey closure' of
$\ch(G)$, or in other words the smallest natural subring of $H^*(G)$
that contains $\ch(G)$ and is closed under transfers.  More formally, 
we make the following:  

\definition Let $G$ be a finite group.  Define $\chbar(G)$ recursively
as the subring of $H^*(G)$ generated by $\ch(G)$ and the image of 
$\chbar(H)$ under the transfer $\Cor^G_H$ for each proper subgroup
$H<G$.  
\medskip

We shall prove: 

\proclaim Theorem 8.1.  Let $G$ be a finite group, with Quillen
category $\AA$, and let $R = \chbar(G)$.  The map $\alpha$ induces
a homeomorphism $\alpha \colon  \colim_\AA X_E(k) \rightarrow V_R(k)$.  

We do not use the injectivity of Quillen's map in proving Theorem~8.1,
so two immediate corollaries of Theorem~8.1 are:  

\proclaim Corollary 8.2. (Quillen) For a finite group $G$, the map
$\alpha \colon  \colim_\AA X_E(k) \rightarrow X_G(k)$ is injective.  
\qed

\proclaim Corollary 8.3.  For a finite group $G$, the inclusion of
$\chbar(G)$ in $H^*(G)$ induces a homeomorphism of varieties.  
\qed

\par\noindent 
{\it Proof of the theorem.}\quad The transfer map $\Cor^G_H$ commutes
with the action of the Steenrod algebra, by either a topological
argument [\ada] or an algebraic one [\evtra]. 
It follows that $R=\chbar(G)$ is a large natural subring of $H^*(G)$.
By Theorem~6.1, it suffices to show that $\CC(R)=\AA$.  

First we show that if $E$, $F$ are elementary abelian subgroups of $G$
such that $E$ is not conjugate to a subgroup of $F$, 
then there is no map in $\CC(R)$ from $E$ to $F$.  Since any map
is the composite of an isomorphism followed by the inclusion of a
subgroup, it suffices to consider the case when $E$ and $F$ have
the same rank.  (Note that if $f$ is any map in $\CC(R)$ that is an
isomorphism of elementary abelian groups, then the inverse of $f$ is
also in $\CC(R)$, so $f$ is an isomorphism in $\CC(R)$.)

Let $N=N_G(F)$ be the normalizer of $F$ in $G$,
let $\theta= \cc F - 1$ be the $(|F|-1)$-dimensional representation of
$F$ given by the difference of the regular representation and the
trivial representation, and let $\rho=\Ind^N_F(\theta)$ be the induced
representation of $N$.  Equivalently, $\rho$ is the regular
representation of $N$ minus the permutation representation on the
cosets $N/F$ of $F$.  Note that $\rho$ is a genuine representation of 
$N$ of dimension $|N|-|N\colon F|$.  Now let $F'$ be any elementary abelian
subgroup of $N$.  The regular representation of $N$ restricts to $F'$
as a sum of $|N\colon F'|$ copies of the regular representation of $F'$.
The number of orbits of $F'$ on the cosets $N/F$, or equivalently the
number of trivial $F'$-summands of the permutation module $\cc N/F$, 
is equal to $|N\colon F'F|$.  It follows that $\Res^N_{F'}(\rho)$, the
restriction to $F'$ of $\rho$, contains the trivial $F'$-module as a
direct summand if and only if $F'F\neq F'$, \ie, if and only if $F$
is not a subgroup of $F'$.  

Let $n$ be the dimension of $\rho$, so that $n= |N|(1-1/|F|)$, and
note that since the restriction to $F$ of $\rho$ does not contain 
the trivial representation, $\Res^N_F(c_n(\rho))$ is non-zero.  
Let $x_F = \Cor^G_N(c_n(\rho))$.  For $E$ an elementary abelian
subgroup of $G$ of the same rank as $F$, 
the Mackey formula affords a calculation of 
$\Res^G_E(x_F)$:  
$$\Res^G_E(x_F)= \sum_{EgN} \Cor^E_{E\cap
gNg^{-1}}c_g^*\Res^N_{g^{-1}Eg\cap N}(c_n(\rho)),$$ 
where for any subgroup $K$ of $G$, $c_g$ is the homomorphism $k\mapsto
g^{-1}kg$, and the sum is over a set of double coset representatives
for $E\backslash G/N$.  The restriction map from $E$ to any subgroup
is surjective in cohomology, and $\Cor^E_{E'}\Res^E_{E'}$ is equal to
multiplication by $|E\colon E'|$.  Hence the
transfer $\Cor^E_{E'}$ is zero for any proper subgroup $E'$ of $E$.
Thus the only non-zero contributions to the above sum come from terms
in which $g^{-1}Eg \leq N$.  On the other hand, we know that the
restriction of $c_n(\rho)$ to an elementary abelian subgroup of $N$ is
non-zero if and only if that subgroup contains $F$.  Since we are
assuming that $E$ has the same rank as $F$, it follows that 
the only non-zero terms will come from $g$ such that $F=g^{-1}Eg$.  
If $F=g^{-1}Eg = h^{-1}Eh$, then $g^{-1}hFh^{-1}g = F$, so $g^{-1}h\in N$, 
and so $EgN= EhN$.  It follows that, for $E$ and $F$ of the same rank,
$\Res^G_E(x_F)= c_g^*(c_n(\rho))$ for any $g$ such that $g^{-1}Eg =
F$, and is zero if there is no such $g$, \ie, if
$E$ and $F$ are not conjugate in $G$.  Since $\Res^N_Fc_n(\rho)$ is
non-zero in $h^*(F)$, it follows that when $E$ and $F$ have the same
rank, there are morphisms from $E$ to $F$ in $\CC(R)$ only if $E$ is
conjugate to $F$.  

It remains to show that the automorphisms of $F$ in the category
$\CC(R)$ are precisely the maps induced by conjugation in $G$.  
Let $C=C_G(F)$ be the centralizer of $F$ in $G$, and suppose that
$|C\colon F|=p^m r$, for some $r$ coprime to $p$.  For any representation
$\lambda$ of $F$, $\Res^C_F\Ind^C_F(\lambda)= p^mr \lambda$.  It
follows that the image of $\ch(C)$ in $\ch(F)$ contains the subring of
$p^m$th powers.  
In $\fp[F]=\ch(F)$, there is a homogeneous element $y_1$ such that the
$\fp\gl(F)$-submodule generated by $y_1$ is free (see [\alp], pp.\
45--46).  The element $y_2=y_1^{p^m}$ also has this property, and 
$y_2=\Res^C_F(y'_2)$ for some $y'_2\in \ch(C)$.  Now let $y'=
y'_2\Res^N_C(c_n(\rho))$, where $\rho$ and $n$ are as in the previous
paragraph.  Then $y'$ is an element of $\ch(C)$ whose restriction to
an elementary abelian subgroup of $C$ is non-zero only if that
subgroup contains $F$.  Moreover, the restriction to $F$ of the
representation $\rho$ is invariant under $\gl(F)$, and so
$y=\Res^C_F(y')$ generates a free $\fp\gl(F)$-submodule of $\ch(F)$.
Now define $z_F\in H^*(G)$ by $z_F = \Cor^G_C(y')$.  The Mackey
formula gives 
$$\Res^G_F(z_F)= \sum_{FgC}\Cor^F_{F\cap
gCg^{-1}}c_g^*\Res^C_{g^{-1}Fg\cap C}(y').$$
Now $\Res^C_{F'}(y') = 0$ unless $F'$ contains $F$, and so 
only those terms for which $g^{-1}Fg=F$ can be non-zero.  Thus 
$$\Res^G_F(z_F)= \sum_{g\in N/C} c_g^*(y),$$ 
where the sum is over cosets of $C=C_G(F)$ in $N= N_G(F)$.  
Since $y$ generates a free $\fp\gl(F)$-summand of $\ch(F)$, 
an automorphism $f$ of $F$ satisfies $f^*\Res^G_F(z_F)=\Res^G_F(z_F)$ 
if and only if $f$ is equal to conjugation by some element of $N$.  
\qed

\remark The first part of this proof is very similar to Quillen's
second proof of this statement, using the Evens norm map
[\pnilp,\evbk,\ben].  The second part is less similar however.  
Our argument is complicated by the weaker technique that we are using
to construct elements, but is simplified by our use of Theorem~6.1 
which means that we do not need to construct as many elements as are
needed in [\pnilp].  

Corollary~8.3 seems to be fairly well-known, although we have been
unable to find it stated in the literature.  Our first proof of
Corollary~8.3 was essentially independent of the rest of this paper,
but used a comparatively recent theorem of Carlson ([\jfc], or
theorem~10.2.1 of [\evbk]):  For $G$ a $p$-group, with centre $Z$, the
radical of $\ker(\Res^G_Z)$ is equal to the radical of the ideal generated
by the images of $\Cor_H^G$, where $H$ ranges over all proper
subgroups of $G$.  

To prove Corollary~8.3 directly, note that if $G_p$ is a Sylow
$p$-subgroup of $G$, the transfer from $H^*(G_p)$ to $H^*(G)$ is
surjective.  Thus it suffices to consider the case when $G$ is a
$p$-group.  Let $G$ be a $p$-group, with centre $Z$.  Using
representations induced from $Z$ up to $G$ it may be shown that the
image of $\ch(G)$ in $H^*(Z)$ contains the subring of $p^m$th powers
for sufficiently large $m$ (as in the proof of Theorem~8.1).  Thus if
$y\in H^*(G)$, there exists $m$ and $x_1\in \ch(G)$ such that 
$y_1= y^{p^m} - x_1$ is in the kernel of $\Res^G_Z$.  By Carlson's
theorem, there exists $n$, subgroups $H(1),\ldots,H(l)$ of $G$ and
$x'_i\in H^*(H(i))$ such that 
$$y_2= y_1^{p^n}= \sum_i\Cor_{H(i)}^G(x'_i).$$
By induction on the order of $G$, there exists $N$ such that, for each
$i$, $(x'_i)^{p^N}\in \chbar(H(i))$.  Noting that for $x$ of degree $2i$, 
$$\Cor^G_H(x^p) = \Cor^G_H(P^ix) = P^i\Cor^G_H(x)= \Cor^G_H(x)^p,$$
it follows that 
$$y_2^{p^N} = \sum_i\Cor_{H(i)}^G((x'_i))^{p^N}=
\sum_i\Cor_{H(i)}^G((x'_i)^{p^N})\in \chbar(G).$$

\beginsection 9. A closure operation.  

\definition
Let $\CC$ be a category of elementary abelian subgroups of a group
$G$.  Define $\CCb$, the closure of $\CC$, to be the smallest
subcategory of $\CC_\reg$ such that:  
\item{1.} $\CC$ is contained in $\CCb$; 

\item{2.} if $f:E_1 \rightarrow E_2$ is in $\CCb$, and $F_i\leq E_i$
with $f(F_1)\leq F_2$, then $f:F_1\rightarrow F_2$ is in $\CCb$; 

\item{3.} if $f:E_1\rightarrow E_2$ is in $\CCb$ and is an isomorphism
of groups, then $f^{-1}: E_2\rightarrow E_1$ is in $\CCb$.  
\medskip

Say that $\CC$ is closed if $\CC = \CCb$.  Note that the categories
$\AA$, $\AA'$, and $\CC(R)$ for any $R\leq H^*(G)$ are closed.  

\proclaim Proposition 9.1.  For any $\CC$ containing the Quillen
category $\AA$, the category $\CCb$ is the 
unique largest subcategory of $\CC_\reg$ such that the natural map 
$$\colim_\CC X_E(k)\rightarrow\colim_\CCb X_E(k)$$ 
is a homeomorphism.  

\proof Let $\DD$ be the subcategory of $\CC_\reg$ whose morphisms
$f:E\rightarrow F$ are those group homomorphisms that make 
the diagram
$$\matrix{X_E(k)&\mapright{f}&X_F(k)\cr 
\mapdown{\ii}&&\mapdown{\ii}\cr
\colim_\CC X_E(k)&\mapright{{\rm Id}}&\colim_\CC X_E(k)\cr}$$
commute.  Then $\DD$ has the property claimed, and it suffices to show
that $\CCb=\DD$.  Note also that $\CCb$ is contained in $\DD$ and that
$\DD$ is closed.  Let $f\colon E_1\rightarrow E_2$ be a morphism in
$\DD$.  Since $f\colon E_1\rightarrow E_2$ is in $\DD$ (resp.\ in
$\CCb$) if and only if $f\colon E_1\rightarrow f(E_1)$ is, it may be
assumed that $f$ is a group isomorphism.  Let $\phi$ be an element of
$X_{E_1}^+(k)$, \ie, an element of $X_{E_1}(k)$ not contained in
$X_F(k)$ for any proper subgroup $F$ of $E_1$.  Since $\gl(E_1)$ acts
freely on $X_{E_1}^+(k)$, it follows that $f\colon E_1\rightarrow E_2$
is uniquely determined by $\psi=f_*(\phi)$.  

By definition of $\DD$, 
$\psi$ and $\phi$ have the same image in $\colim_\CC X_E(k)$.  
Since $\CC$ is skeletally finite (because it contains $\AA$), 
there are chains $(F_0,\ldots,F_m)$ of objects of $\CC$ and
$(f_1,\ldots,f_m)$ of morphisms in $\CC$, where 
$$f_i\colon F_{i-1+\epsilon(i)}\rightarrow F_{i-\epsilon(i)}$$
for some
$\epsilon(i)\in\{0,1\}$, and $(\psi_0,\ldots,\psi_m)$, $\psi_i\in
X_{F_i}(k)$, with 
$$F_0=E_1,\qquad F_m=E_2,\qquad \psi_0=\phi,\qquad \psi_m=\psi, 
\qquad {f_i}_*(\psi_{i-1+\epsilon(i)})= \psi_{i-\epsilon(i)}.$$
Let $F'_i$ be the unique subgroup of $F_i$ such that $\psi_i\in
X_{F'_i}^+(k)$.  Then $F'_i$ has the same rank as $E_1$ and $f_i$
restricts to an isomorphism $f'_i$ from $F'_{i-1+\epsilon(i)}$ to
$F'_{i-\epsilon(i)}$.  Letting $\delta(i)= 1-2\epsilon(i)$,
$f_i^{\prime\delta(i)}$ is a morphism in $\CCb$ from $F'_{i-1}$ to
$F'_i$, and the composite 
$$f'= f_m^{\prime\delta(m)}\circ\cdots\circ f_1^{\prime\delta(1)}$$ 
is a morphism in $\CCb$ from $E_1$ to $E_2$ such that
$f'(\phi)=\psi$.  Hence $f'=f$, and $f$ is a morphism in $\CCb$ as
claimed.  
\qed

For any category $\CC$ of elementary abelian subgroups of a group $G$,
one may define a subring $R(\CC)$ of $H^*(G)$ as the inverse image of 
$\lim_\CC H^*(E)$.  This subring is large and is natural 
because $\lim_\CC H^*(E)$ is.  

\proclaim Proposition 9.2.  For $\CC$ any category of elementary
abelian subgroups of $G$ containing the Quillen category $\AA$, 
$\CC(R(\CC))=\CCb$.  

\par\noindent
\proof Clearly, $\CC(R(\CC))$ contains $\CCb$, and 
is closed.  Hence it suffices to show that the induced map of
varieties is a homeomorphism.  Quillen showed that the
map from $H^*(G)$ to $\lim_\AA H^*(E)$ contains the subring of
$p^n$th powers for some $n$ (in fact this is equivalent to the
injectivity of the map $\colim_\AA X_E(k)\rightarrow X_G(k)$) [\qui]. 
Let $S= \lim_\CC H^*(E)$, and note that if $x$ is any element of $S$,
the $p^n$th power of $x$ is in the image of $R=R(\CC)$.  
It follows that the
map $V_S(k) \rightarrow V_R(k)$ is a homeomorphism, and it suffices to
show that the natural map 
$$\colim_\CC X_E(k) \rightarrow V_S(k)$$ 
is a homeomorphism.  But this is a special case of lemma~8.11 of
[\qui].  \qed

The proposition shows that there is a sort of \lq Galois
correspondence' between large natural subrings of $H^*(G)$ and
categories of elementary abelian subgroups of $G$.  

\beginsection 10. Some other categories.  

For each $n\geq 0$, define a category $\AA^{(n)}$ of elementary 
abelian subgroups of a group $G$ by declaring that the morphism 
$f \colon E\hookrightarrow F$ is in $\AA^{(n)}$ if and only if
for all $e_1,\ldots,e_n\in E$, there exists $g\in G$ such that 
$f(e_i)= g^{-1}eg$.  
\medskip

Note that $\AA^{(0)}$ is the category $\CC_\reg$
of Section~6, and $\AA^{(1)}$ is the category $\AA'$.  For each~$n$, 
$\AA^{(n)}\supseteq \AA^{(n+1)}$, and when $n$ is greater than or
equal to the $p$-rank of $G$, $\AA^{(n)}$ is equal to Quillen's
category $\AA$.  This suggests that $\AA^{(\infty)}$ should be defined
to be $\AA$.  Each of the categories $\AA^{(n)}$ is closed in the
sense of section~9, and the subrings $R^{(n)}=R(\AA^{(n)})$ form a
natural filtration of $H^*(G)=R^{(\infty)}$.  

The categories $\AA^{(n)}(G)$ are related to the generalized
characters of $G$ due to Hopkins, Kuhn and Ravenel [\hkr] in the same 
way that the category $\AA'$ is related to ordinary characters.  
It seems possible that there should be a description of the variety
for the subring of elements of $H^*(G)$ coming from $E^0(BG)$, where
$E$ is a generalized cohomology theory to which Hopkins-Kuhn-Ravenel's
work applies, in the same way that Chern classes are elements of
$H^*(G)$ coming from $K^0(BG)$.  We shall not make a precise 
conjecture, but shall give examples to show that the
categories $\AA^{(n)}$ can be distinct from each other.

\proclaim Proposition 10.1.  For each $n\geq 0$ and each prime $p$,
there is a $p$-group $G$ for which $\AA^{(n)}(G)\neq 
\AA^{(n+1)}(G)$.  

\proof For $n=0$, the cyclic group of order $p$ (for $p$ odd), or 
the elementary abelian group of order four (for $p=2$), will suffice.  
Hence we may assume that $n>0$.  
Let $C$ be a cyclic group of order $p$, let $E$ be
a faithful $\fp C$-module of $\fp$-dimension $n+1$, and let $c\in 
\gl(E)$ represent the action on $E$ of a generator for $C$.  Now let $Z$ 
be a vector space over $\fp$ with basis $\{z_M\}$ indexed by the maximal
$\fp$-subspaces of $E$, so that $Z$ has dimension $(p^{n+1}-1)/(p-1)$.  

For each maximal subspace $M$ of $E$, pick a linear map 
$\psi_M\colon E\rightarrow Z$, with kernel $M$ and image generated by
$z_M$.  For each $M$, define $b_M\in \gl(E\oplus Z)$ by the equation  
$$b_M(e,z) = (c(e),z+\psi_M(e)).$$ 
Let $A$ be the subgroup of $\gl(E\oplus Z)$ generated by the $b_M$,
and let $G$ be the semidirect product $(E\oplus Z)\colon A$.  

The subgroup $Z$ is left invariant by $A$, so is central in $G$.  
Let $\phi$ be the homomorphism sending $A\leq \gl(E\oplus Z)$ to 
$\gl((E\oplus Z)/Z)\cong \gl(E)$, and let $B=\ker(\phi)\leq A$.  
Note that elements of $B$ act trivially on $Z$ and on $(E\oplus Z)/Z$,
and so $B$ may be identified with a subgroup of the elementary abelian
$p$-group $\Hom(E,Z)$.  

We claim that the automorphism $c$ of $E$ is a morphism in
$\AA^{(n)}(G)$, but is not a morphism in $\AA^{(n+1)}(G)$.  If $M$ is
any rank-$n$ subgroup of $E$, then the element $b_M\in G$ acts on $M$
in the same way as $c$.  On the other hand, if $c$ were a morphism in
$\AA^{(n+1)}(G)$, there would have to be an element $d$ of $G$, acting
on $E\oplus Z$ as $d(e,z)=(c(e),z)$.  But then, for any $M$,
$d'=d^{-1}b_M$ would be an element of $(E\oplus Z):B$ acting as
$d'(e,z)=(e,z+\psi_M(e))$.  To complete the proof, it suffices to show
that there can be no such element $d'$.  

Let $R$ be the image of $\fp C$ in the ring $\End(E)$.  Since $\fp C$
is a commutative local ring, it follows that $R$ is too.  In
particular, the non-units in $R$ form an ideal.  Fix $M$, a maximal
subgroup of $E$.  The group $B\leq A
\leq \gl(E)$ is generated as a normal subgroup by the elements
$b_M^p$, and $b_M^{-1}b_N$, where $N$ ranges over all other maximal
subgroups of $E$.  The action of these generators on $E\oplus Z$ is
given by:  
$$b_M^p(e,z)= \left(c^p(e),z+\sum_{i=0}^{p-1} \psi_M(c^i(e))\right) 
= (e,\psi_M(\bar r e)+z),$$
where $\bar r$ is the image of $\bar c= \sum_{i=0}^{p-1} c^i$ in
$R=\End(E)$, and 
$$b_M^{-1}b_N(e,z)= (e,z+\psi_N(e) - \psi_M(e)).$$
We therefore have to show that the element $d'$ described above does
not lie in the subgroup of $\gl(E\oplus Z)$ generated by 
$$(e,z)\mapsto (e,z+\psi_M(\bar r e))\qquad 
\hbox{and}\qquad (e,z)\mapsto (e,z+\psi_N(c^i e) - \psi_M(c^ie)),$$
for all $N\neq M$ and $0\leq i\leq p-1$.  

If $d'\colon (e,z)\mapsto (e,z+\psi_M(e))$ is in the subgroup $B$,
there are $\mu,\lambda_N\in R\leq \End(E)$ such that 
$$\Im \lambda_N \subseteq \ker(\psi_N)= N\quad \hbox{for all $N\neq
M$, and}$$
$$\Im \left(1 -\mu \bar r + \sum_{N\neq M} \lambda_N\right)
\subseteq \ker(\psi_M)= M.$$
From the first family of equations, it follows that each $\lambda_N$
is a non-unit in $R$, and of course $\bar r$ is a non-unit.  But then 
$1-\mu \bar r + \sum \lambda_N$ must be a unit in $R$, contradicting
the second equation.  
\qed

\remark We exhibit $p$-group examples above because $p$-groups 
control the behaviour of mod-$p$ cohomology, and hence $p$-group
examples tend to be harder to find.  Similar examples can be
constructed for other cyclic groups $C$.  One interesting example is
the case when $C=\gl_1(\fq)$ for $q=p^{n+1}>2$ and $E$ is the additive
group of $\fq$.  In this case $\fp C$ is not a local ring, but its
image in $\End(E)$ is isomorphic to $\fq$, and $\bar r = 0$, so the
argument used above still applies.

\beginsection References.  

\frenchspacing 

\def\book#1/#2/#3/#4/ {\item{[#1]} #2, {\it #3}, #4.
\par\smallskip}
\def\paper#1/#2/#3/#4/#5/ {\item{[#1]} #2, {\it #3} {\bf #4}, #5.
\par\smallskip}

\book\ada/J. F. Adams/Infinite loop spaces/Annals of Mathematics
Studies vol. 90, Princeton University Press, Princeton (1978)/ 

\book\alp/J. L. Alperin/Local representation theory/ Cambridge
Studies in Advanced Mathematics vol. 11, Cambridge University Press,
Cambridge (1986)/

\paper\ati/M. F. Atiyah, Characters and cohomology of finite
groups/Publ. Math. IHES/9/(1961), 23--64/ 

\book\atmc/M. F. Atiyah and I. G. Macdonald/Introduction to
commutative algebra/Addison-Wesley, Reading, Mass. (1969)/

\book\ben/D. J. Benson/Representations and cohomology II: cohomology
of groups and modules/Cambridge Studies in Advanced Mathematics
vol. 31, Cambridge University Press, Cambridge (1991)/ 

\book\beninv/D. J. Benson/Polynomial invariants of finite
groups/London Math. Soc. Lecture Notes vol. 190, Cambridge University
Press, Cambridge (1993)/  

\paper\jfc/J. F. Carlson, Varieties and transfer/J. Pure
Appl. Alg./44/(1987) 99--105/

\paper\dww/W. G. Dwyer and C. W. Wilkerson, Homotopy fixed point
methods for Lie groups and finite loop spaces/Annals of 
Math./139/(1994) 395--442/

\paper\ev/L. Evens, The cohomology ring of a finite
group/Trans. Amer. Math. Soc./101/(1961), 224--239/ 

\paper\evtra/L. Evens, Steenrod operations and
transfer/Proc. Amer. Math. Soc./19/(1968), 1387--1388/

\book\evbk/L. Evens/The cohomology of groups/Oxford Mathematical
Monographs, Clarendon Press, Oxford (1991)/  

\paper\hls/H.-W. Henn, J. Lannes and L. Schwartz, The categories of
unstable modules and unstable algebras modulo nilpotent
objects/Amer. J. Math/115/(1993), 1053--1106/  

\book\lam/S. P. Lam/Unstable algebras over the Steenrod
algebra/Lecture Notes in Mathematics vol. 1051, Springer, New York
(1984)/  

\item{[\boaz]} B. Moselle, Unpublished essay, University of Cambridge
(1989). 
\smallskip

\paper\qui/D. Quillen, The spectrum of an equivariant cohomology ring
I and II/Ann. of Math./94/(1971), 549--572 and 573--602/ 

\paper\pnilp/D. Quillen, A criterion for $p$-nilpotence/J. Pure
Appl. Alg./1/(1971) 361--372/

\paper\rec/D. Rector, Noetherian cohomology rings and finite loop
spaces with torsion/J. Pure Appl. Alg./32/(1984) 191--217/

\paper\ser/J.-P. Serre, Sur la dimension cohomologique des groupes
profinis/Topology/3/(1965) 413--420/  

\book\serbk/J.-P. Serre/Linear representations of finite
groups/Graduate Texts in Mathematics vol. 42, Springer, New York
(1977)/  

\book\tho/C. B. Thomas/Characteristic classes and the cohomology of
finite groups/Cam\-bridge Studies in Advanced Mathematics vol. 9,
Cambridge University Press, Cambridge (1986)/  

\paper\ven/B. B. Venkov, Cohomology algebras for some classifying
spaces/Dokl. Akad. Nauk SSSR/127/(1959), 943--944/

\item{[\wilk]} C. W. Wilkerson, A primer on the Dickson invariants, in
Proc. Northwestern Homotopy Theory Conference, Contemp. Math vol. 19
(1983), 421--434.  
\smallskip 


\vskip .5in 
$$\hbox{
\vbox{
\hbox{David J. Green} 
\hbox{Inst. f. Experimentelle Math.,}
\hbox{Ellernstra\ss e 29,}
\hbox{D-45326 Essen}
\hbox{Germany}
\hbox{\phantom{United Kingdom}}
}\qquad\vbox{
\hbox{Ian J. Leary}
\hbox{Faculty of Math. Studies,}
\hbox{Univ. of Southampton,}
\hbox{Southampton,}
\hbox{SO17 1BJ} 
\hbox{United Kingdom}
}}$$

\end